\documentclass[12pt]{article}%
\usepackage{amsfonts}
\usepackage{sw20bams}
\usepackage{amsmath}
\usepackage{amssymb}
\usepackage{graphicx}%
\providecommand{\U}[1]{\protect\rule{.1in}{.1in}}
\begin{document}

\title{Lognormality and Triangles of Unit Area}
\author{Steven Finch}
\date{December 18, 2014}
\maketitle

\begin{abstract}
To generate a triangle of unit perimeter, break a stick of length $1$ in two
places at random, with the condition that triangle inequalities are satisfied.
\ Is there a similarly natural method for generating triangles of unit area?
\ Study of a\ product (rather than a sum) of random lengths is facilitated by
closure of multiplication within the lognormal family of distributions. \ Our
(necessarily incomplete)\ answers to the question each draw upon this property.

\end{abstract}

\footnotetext{Copyright \copyright \ 2014 by Steven R. Finch. All rights
reserved.}We denote triangle sides by $a$, $b$, $c$ and opposite angles by
$\alpha$, $\beta$, $\gamma$. \ The following formulas for area:%
\[
\frac{1}{2}a\,b\sin(\gamma)=\frac{1}{4}\sqrt{(a+b+c)(-a+b+c)(a-b+c)(a+b-c)}
\]
will be used (out of many possibilities \cite{Bk-LogNrm}). \ Also, linear
area-preserving transformations $\mathbb{R}^{2}\rightarrow\mathbb{R}^{2}$
constitute the group $SL_{2}(\mathbb{R)}$ of all real $2\times2$ matrices with
determinant $1$. Every such matrix possesses a unique representation as
\cite{Cn-LogNrm}%
\[
\left(
\begin{array}
[c]{ll}%
\cos(\theta) & -\sin(\theta)\\
\sin(\theta) & \cos(\theta)
\end{array}
\right)  \left(
\begin{array}
[c]{ll}%
r & 0\\
0 & 1/r
\end{array}
\right)  \left(
\begin{array}
[c]{ll}%
1 & s\\
0 & 1
\end{array}
\right)
\]
where $0\leq\theta<2\pi$, $r>0$ and $s\in\mathbb{R}$. \ We will use only the
center matrix $M_{r}$ (with $r$ and $1/r$ on the diagonal) in this paper.
\ Incorporating the right matrix (with $s$ in the off-diagonal corner) is left
for future research.

\section{Right Triangles}

Assuming $\gamma=\pi/2$, the sides $a$, $b$ must satisfy the constraint
$a\,b=2$. \ Given $a$, it is clear that $b$, $c$ are completely determined
and, given $\alpha$, we have $\beta=\pi/2-\alpha$. \ From%
\[
\ln(a)+\ln(b)=\ln(2)
\]
it becomes reasonable to adopt the model $\ln(a)\sim\,$Normal$\left(  \frac
{1}{2}\ln(2),\sigma^{2}\right)  $. By definition, this is equivalent to
$a\sim\,$Lognormal$\left(  \frac{1}{2}\ln(2),\sigma^{2}\right)  $. \ The
density function for $a$ is%
\[%
\begin{array}
[c]{lll}%
\dfrac{1}{\sqrt{2\pi}\sigma}\exp\left[  -\dfrac{1}{2\sigma^{2}}\left(
\ln(a)-\dfrac{1}{2}\ln(2)\right)  ^{2}\right]  \dfrac{1}{a}, &  & a>0
\end{array}
\]
with moments
\[%
\begin{array}
[c]{lll}%
\operatorname*{E}\left(  a\right)  =\sqrt{2e^{\sigma^{2}}}, &  &
\operatorname*{E}\left(  a^{2}\right)  =2e^{2\sigma^{2}}.
\end{array}
\]
For convenience, we set $\sigma=1$ always. \ From the Law of Sines,%
\[
\sin(\alpha)=\frac{a}{b}\sin(\beta)=\frac{a}{2/a}\sin\left(  \frac{\pi}%
{2}-\alpha\right)  =\frac{a^{2}}{2}\cos(\alpha)
\]
hence%
\[%
\begin{array}
[c]{lll}%
2\tan(\alpha)=a^{2}, &  & 2\sec(\alpha)^{2}d\alpha=2a\,da
\end{array}
\]
hence%
\[
\dfrac{da}{a}=\dfrac{d\alpha}{a^{2}\cos(\alpha)^{2}}
\]
thus the density function for $\alpha$ is%
\begin{align*}
&  \dfrac{1}{\sqrt{2\pi}}\exp\left[  -\dfrac{1}{2}\left(  \ln\left(
\sqrt{2\tan(\alpha)}\right)  -\dfrac{1}{2}\ln(2)\right)  ^{2}\right]
\dfrac{1}{2\tan(\alpha)\cos(\alpha)^{2}}\\
&  =\dfrac{1}{\sqrt{2\pi}}\exp\left[  -\dfrac{1}{8}\ln\left(  \tan
(\alpha)\right)  ^{2}\right]  \dfrac{1}{2\sin(\alpha)\cos(\alpha)}%
\end{align*}
for $0<\alpha<\pi/2$. \ Let $x=\ln\left(  \tan(\alpha)\right)  /2$, then
$\alpha=\arctan\left(  e^{2x}\right)  $ and
\[
dx=\frac{\sec(\alpha)^{2}}{2\tan(\alpha)}d\alpha=\frac{d\alpha}{2\sin
(\alpha)\cos(\alpha)}
\]
therefore the moments become%
\begin{align*}
\operatorname*{E}\left(  \alpha\right)   &  =\dfrac{1}{\sqrt{2\pi}}%
{\displaystyle\int\limits_{-\infty}^{\infty}}
\arctan\left(  e^{2x}\right)  e^{-x^{2}/2}dx\\
&  =\frac{\pi}{4}\\
&  =0.7853981633974483096156608...,
\end{align*}%
\begin{align*}
\operatorname*{E}\left(  \alpha^{2}\right)   &  =\dfrac{1}{\sqrt{2\pi}}%
{\displaystyle\int\limits_{-\infty}^{\infty}}
\arctan\left(  e^{2x}\right)  ^{2}e^{-x^{2}/2}dx\\
&  =0.9012156209647814268211368....
\end{align*}
No closed-form expression for the mean square of $\alpha$ is known.%
\begin{figure}[ptb]%
\centering
\includegraphics[
height=2.8885in,
width=6.1307in
]%
{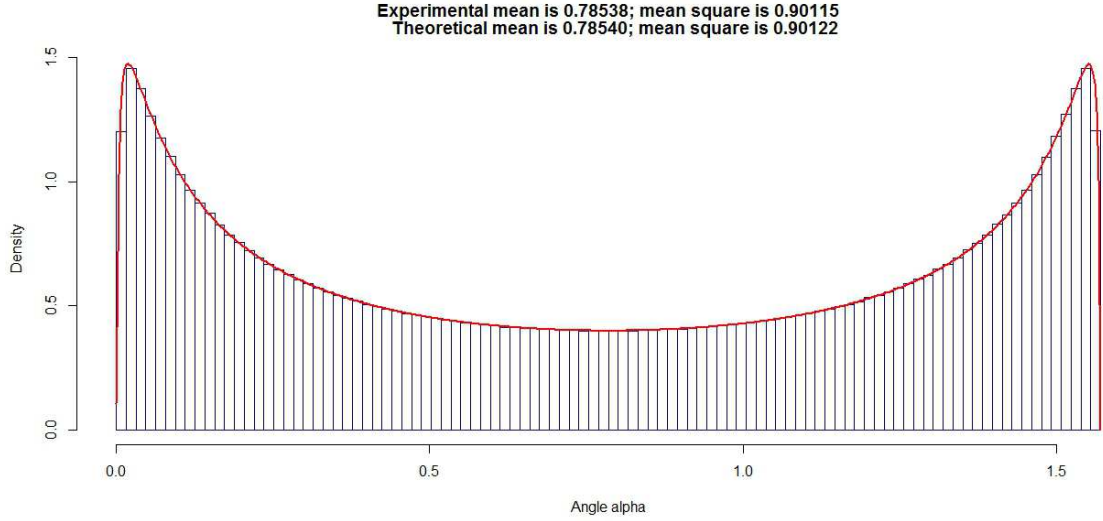}%
\caption{Density function for angle $\alpha$ in Section 1; there are zeroes at
$0$ and $\pi/2$ but maximum points at $\varepsilon$ and $\pi/2-\varepsilon$,
where $\varepsilon\approx0.018363$.}%
\end{figure}

\section{Isosceles Triangles}

The equilateral triangle $T$ with vertices
\[%
\begin{array}
[c]{ccccc}%
A=\left(  \dfrac{1}{3^{1/4}},-\dfrac{1}{3^{3/4}}\right)  , &  & B=\left(
-\dfrac{1}{3^{1/4}},-\dfrac{1}{3^{3/4}}\right)  , &  & C=\left(  0,\dfrac
{2}{3^{3/4}}\right)
\end{array}
\]
is centered at $(0,0)$ and has area $1$. \ Images of $T$ under $M_{r}$ for
$r>0$ encompass all unit area triangles satisfying the constraint $a=b$. \ It
is reasonable to adopt the model $r\sim\,$Lognormal$\left(  -\frac{1}%
{2},1\right)  $ so that $\operatorname{E}(M_{r})$ is the identity matrix.
\ Clearly%
\[
a=\left\vert M_{r}(C-B)\right\vert =\sqrt{\frac{r^{2}}{\sqrt{3}}+\frac
{\sqrt{3}}{r^{2}}},
\]%
\[
b=\left\vert M_{r}(C-A)\right\vert =\sqrt{\frac{r^{2}}{\sqrt{3}}+\frac
{\sqrt{3}}{r^{2}}},
\]%
\[
c=\left\vert M_{r}(B-A)\right\vert =\sqrt{\frac{4r^{2}}{\sqrt{3}}}
\]
are the sides and consequently $c\sim\,$Lognormal$\left(  -\frac{1}{2}%
+\ln\left(  \tfrac{2}{3^{1/4}}\right)  ,1\right)  $,%
\[%
\begin{array}
[c]{lll}%
\operatorname*{E}\left(  c\right)  =\dfrac{2}{3^{1/4}}, &  & \operatorname*{E}%
\left(  c^{2}\right)  =\dfrac{4e}{\sqrt{3}}.
\end{array}
\]
In Section 4, we calculate the density function for $a$ to be%
\[
\frac{a}{\sqrt{2\pi}}\frac{\exp\left\{  -\dfrac{1}{8}\left[  \ln\left(
\dfrac{\sqrt{3}}{2}\left(  a^{2}-\sqrt{a^{4}-4}\right)  \right)  +1\right]
^{2}\right\}  +\exp\left\{  -\dfrac{1}{8}\left[  \ln\left(  \dfrac{\sqrt{3}%
}{2}\left(  a^{2}+\sqrt{a^{4}-4}\right)  \right)  +1\right]  ^{2}\right\}
}{\sqrt{a^{4}-4}}
\]
when $a>\sqrt{2}$. \ Thus the moments of $a$ are%
\[
\operatorname*{E}\left(  a\right)  =3.9753634096801809039980060...,
\]%
\begin{align*}
\operatorname*{E}\left(  a^{2}\right)   &  =\sqrt{3}\left(  \frac{1}{3}%
+e^{2}\right)  e\\
&  =36.3585711936558006990230225....
\end{align*}
No closed-form expression for the mean of $a$ is known. \ We have not
attempted to find any densities or moments for corresponding angles.

\section{Arbitrary Triangles}

A disappointing aspect of the preceding two sections is the lack of any
bivariate densities or cross-correlations. \ We shall now remedy this
situation, but a disadvantage of our third model is its artificiality.

The sides $a$, $b$, $c$ of a triangle $\Delta$ satisfy%
\[
16=(a+b+c)(-a+b+c)(a-b+c)(a+b-c)
\]
if and only if $\Delta$ has area $1$. \ The set of all solutions $(a,b,c)$ of
this quartic equation is a noncompact surface $\Sigma$ in the positive octant
of $\mathbb{R}^{3}$. \ The orthogonal projection of $\Sigma$ into the plane
$c=0$ is the region bounded (away from the origin) by the hyperbola%
\[
a\,b=2.
\]
For example, the interior point $(a,b)=(2\cdot3^{-1/4},2\cdot3^{-1/4})$ of the
region corresponds to the equilateral triangle ($c=2\cdot3^{-1/4}$). \ But the
correspondence is not one-to-one: the point \textit{also} corresponds to an
isosceles triangle with apex angle $\gamma=2\pi/3$\ ($c=2\cdot3^{1/4}$).
\ More generally, we have
\[
c=\sqrt{a^{2}+b^{2}\pm2\sqrt{a^{2}b^{2}-4}}
\]
and the two possible values of $c$ are assigned equal probability in our
model. The only exceptions lie on the hyperbola itself: any boundary point of
the region corresponds uniquely to a right triangle ($c=\sqrt{a^{2}+b^{2}} $).%
\begin{figure}[ptb]%
\centering
\includegraphics[
trim=0.000000in -0.018867in 0.000000in 0.018867in,
height=2.9788in,
width=6.3179in
]%
{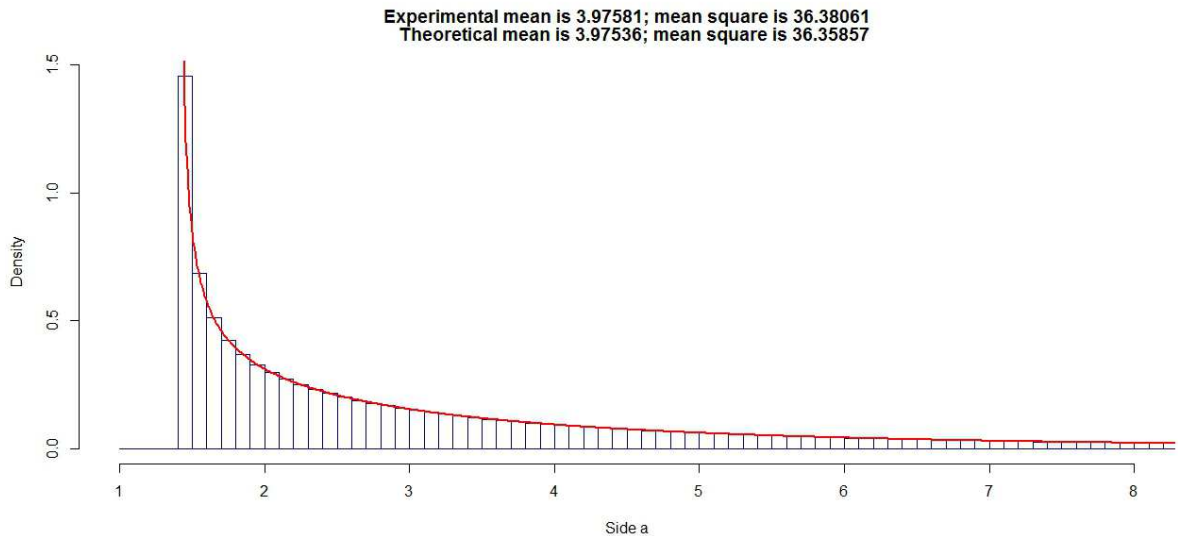}%
\caption{Density function for side $a$ in Section 2; the tail is
asymptotically $O\left(  \exp\left(  -\frac{1}{2}\ln(a)^{2}\right)
/a^{\delta}\right)  $ as $a\rightarrow\infty$, where $\delta=1/2-\ln
(3)/4\approx0.225347.$}%
\end{figure}

Random sampling of $a\,b\geq2$ is performed as follows. \ Generate $u$, $v$
independently according to Normal$\left(  0,1\right)  $. \ If $u+v<0$, then
reflect the point $(u,v)$ in the $xy$-plane across the diagonal line $x+y=0$;
otherwise do nothing. \ [The reflection is achieved via overwriting $(u,v)$ by
$(-v,-u)$.] \ Now translate $(u,v)$ by adding $\frac{1}{2}\ln(2)$ to each
component. The density of $(u,v)$ is thus a folded bivariate normal:%
\[
\frac{1}{\pi}\exp\left\{  -\frac{1}{2}\left[  \left(  u-\frac{1}{2}%
\ln(2)\right)  ^{2}+\left(  v-\frac{1}{2}\ln(2)\right)  ^{2}\right]  \right\}
\]
for $u+v\geq\ln(2)$. \ Let $a=e^{u}$ and $b=e^{v}$, then the bivariate density
of $(a,b)$ is%
\[
\frac{1}{\pi}\exp\left\{  -\frac{1}{2}\left[  \left(  \ln(a)-\frac{1}{2}%
\ln(2)\right)  ^{2}+\left(  \ln(b)-\frac{1}{2}\ln(2)\right)  ^{2}\right]
\right\}  \frac{1}{a\,b}%
\]
for $a\,b\geq2$. \ Integrating on $b$ from $2/a$ to infinity, we obtain the
univariate density of $a$:%
\[
\dfrac{1}{\sqrt{2\pi}}\exp\left[  -\dfrac{1}{2}\left(  \ln(a)-\dfrac{1}{2}%
\ln(2)\right)  ^{2}\right]  \operatorname*{erfc}\left[  -\frac{1}{\sqrt{2}%
}\left(  \ln(a)-\dfrac{1}{2}\ln(2)\right)  \right]  \dfrac{1}{a}%
\]
for $a>0$, where $\operatorname*{erfc}$ is the complementary error function.
\ These results give rise to
\[
\operatorname*{E}\left(  a\right)  =3.5452643891219526811143352...,
\]%
\[
\operatorname*{E}\left(  a^{2}\right)  =27.2316390652988719486867211...,
\]%
\[
\operatorname*{E}\left(  a\,b\right)  =10.0179601615245669326196491...
\]
and, in particular, the cross-correlation coefficient between $a$ and $b$ is
$\approx-0.174$. \ Integrating the average of the two $c$-values%
\[
\dfrac{1}{2}\left(  \sqrt{a^{2}+b^{2}+2\sqrt{a^{2}b^{2}-4}}+\sqrt{a^{2}%
+b^{2}-2\sqrt{a^{2}b^{2}-4}}\right)
\]
gives $\operatorname*{E}\left(  c\right)  \approx5.483$, which unfortunately
demonstrates how artificial this model is. \ We wish ideally for
$\operatorname*{E}\left(  c\right)  $ to be equal to both $\operatorname*{E}%
\left(  a\right)  $ and $\operatorname*{E}\left(  b\right)  $, since there is
no reason for one side to be preferred over the other two. \ An alternative
approach would involve the orthogonal projection of $\Sigma$ into the plane
$c-a=0$ (rather than $c=0$), which may provide the desired symmetry.

Our goal of generating unit-area triangles analogously to unit-perimeter
triangles -- \textquotedblleft throwing paint\textquotedblright\ rather than
breaking a stick -- remains elusive \cite{Fi1-LogNrm}. \ A final comment
concerning Section 2 \&\ Section 4 in \cite{Fi2-LogNrm} (based partly on
\cite{ES-LogNrm}) is in order. \ If the lengths of the three pieces (from
randomly breaking the stick twice) instead are $a^{2}$, $b^{2}$, $c^{2}$ and
all triangle inequalities are satisfied, then area has density function $96x$
over $[0,\sqrt{3}/12]$. \ An even simpler outcome emerges if the lengths of
the two pieces (from breaking the stick just once) are $a^{2}$, $b^{2}$ and
angle $\gamma$ is taken to be Uniform$[0,\pi]$. \ Area in this case is
distributed according to Uniform$[0,1/4]$. \ We wonder whether some elementary
modifications of either case might lead to insight necessary to answer our
question. \
\begin{figure}[ptb]%
\centering
\includegraphics[
height=4.907in,
width=4.0923in
]%
{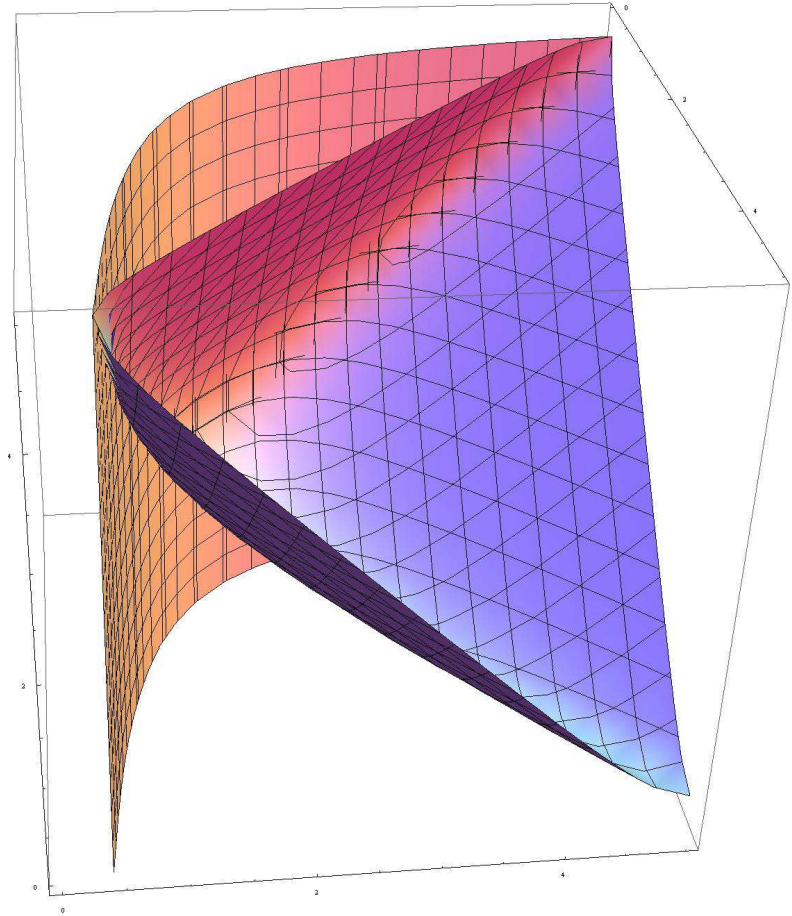}%
\caption{Surface $\Sigma$ in foreground and vertical cylinder (determined by
$a\,b=2$) in background; intersection captures all right triangles.}%
\end{figure}

\section{Details of Calculation}

We seek the distribution of the sum of a random variable $x\sim\,$%
Lognormal$\left(  \mu,\sigma^{2}\right)  $ and its reciprocal $1/x\sim
\,$Lognormal$\left(  -\mu,\sigma^{2}\right)  $. \ Solving the equation%
\[
x+\frac{1}{x}=y
\]
for $x$ in terms of $y$, we obtain two positive values%
\[%
\begin{array}
[c]{ccc}%
x_{-}=\dfrac{1}{2}\left(  y-\sqrt{y^{2}-4}\right)  <1, &  & x_{+}=\dfrac{1}%
{2}\left(  y+\sqrt{y^{2}-4}\right)  >1
\end{array}
\]
for $y>2$. \ Because%
\[%
\begin{array}
[c]{ccccc}%
\dfrac{d}{dx}\left(  x+\dfrac{1}{x}\right)  =1-\dfrac{1}{x^{2}}, &  &
1-\dfrac{1}{x_{-}^{2}}<0, &  & 1-\dfrac{1}{x_{+}^{2}}>0
\end{array}
\]
and%
\begin{align*}
1-\frac{1}{x_{\mp}^{2}}  &  =1-\dfrac{4}{\left(  y\mp\sqrt{y^{2}-4}\right)
^{2}}\\
&  =\dfrac{\left(  y\mp\sqrt{y^{2}-4}\right)  ^{2}-4}{\left(  y\mp\sqrt
{y^{2}-4}\right)  ^{2}}\\
&  =\dfrac{y^{2}\mp2y\sqrt{y^{2}-4}+\left(  y^{2}-4\right)  -4}{\left(
y\mp\sqrt{y^{2}-4}\right)  ^{2}}\\
&  =\dfrac{2y^{2}\mp2y\sqrt{y^{2}-4}-8}{\left(  y\mp\sqrt{y^{2}-4}\right)
^{2}}\\
&  =\dfrac{\mp2\left(  y\mp\sqrt{y^{2}-4}\right)  \sqrt{y^{2}-4}}{\left(
y\mp\sqrt{y^{2}-4}\right)  ^{2}}\\
&  =\dfrac{\mp2\sqrt{y^{2}-4}}{y\mp\sqrt{y^{2}-4}}=\frac{\mp\sqrt{y^{2}-4}%
}{x_{\mp}},
\end{align*}
it follows that the density function for $y$:%
\[
-\dfrac{1}{\sqrt{2\pi}\sigma}\frac{\exp\left\{  -\dfrac{1}{2\sigma^{2}}\left[
\ln\left(  x_{-}\right)  -\mu\right]  ^{2}\right\}  \dfrac{1}{x_{-}}}%
{1-\dfrac{1}{x_{-}^{2}}}+\dfrac{1}{\sqrt{2\pi}\sigma}\frac{\exp\left\{
-\dfrac{1}{2\sigma^{2}}\left[  \ln\left(  x_{+}\right)  -\mu\right]
^{2}\right\}  \dfrac{1}{x_{+}}}{1-\dfrac{1}{x_{+}^{2}}}
\]
simplifies to%
\[
\dfrac{1}{\sqrt{2\pi}\sigma}\frac{\exp\left\{  -\dfrac{1}{2\sigma^{2}}\left[
\ln\left(  \dfrac{1}{2}\left(  y-\sqrt{y^{2}-4}\right)  \right)  -\mu\right]
^{2}\right\}  +\exp\left\{  -\dfrac{1}{2\sigma^{2}}\left[  \ln\left(
\dfrac{1}{2}\left(  y+\sqrt{y^{2}-4}\right)  \right)  -\mu\right]
^{2}\right\}  }{\sqrt{y^{2}-4}}.
\]
Consider now $z=\sqrt{y}$, for which $y=z^{2}$ and%
\[
\frac{d}{dy}\sqrt{y}=\frac{1}{2\sqrt{y}}=\frac{1}{2z};
\]
the density function for $z$ is%
\[
\sqrt{\frac{2}{\pi}}\frac{z}{\sigma}\frac{\exp\left\{  -\dfrac{1}{2\sigma^{2}%
}\left[  \ln\left(  \dfrac{1}{2}\left(  z^{2}-\sqrt{z^{4}-4}\right)  \right)
-\mu\right]  ^{2}\right\}  +\exp\left\{  -\dfrac{1}{2\sigma^{2}}\left[
\ln\left(  \dfrac{1}{2}\left(  z^{2}+\sqrt{z^{4}-4}\right)  \right)
-\mu\right]  ^{2}\right\}  }{\sqrt{z^{4}-4}}
\]
when $z>\sqrt{2}$. \ As an example, if $\mu=-1/2$ and $\sigma=1$, then%
\[
\operatorname*{E}\left(  z\right)  =1.8366252372930300853898532...,
\]%
\begin{align*}
\operatorname*{E}\left(  z^{2}\right)   &  =1+e\\
&  =3.7182818284590452353602874....
\end{align*}
No closed-form expression for the mean of $z=\sqrt{x+1/x}$ is known. \ An
attractive integral representation%
\begin{align*}
\operatorname*{E}\left(  z\right)   &  =\dfrac{1}{\sqrt{2\pi}}%
{\displaystyle\int\limits_{0}^{\infty}}
\sqrt{x+\frac{1}{x}}\exp\left[  -\dfrac{1}{2}\left(  \ln(x)+\dfrac{1}%
{2}\right)  ^{2}\right]  \dfrac{1}{x}dx\\
&  =\dfrac{1}{\sqrt{\pi}}%
{\displaystyle\int\limits_{-\infty}^{\infty}}
\sqrt{\cosh\left(  u-\frac{1}{2}\right)  }\exp\left(  -\frac{1}{2}%
u^{2}\right)  du
\end{align*}
(found via $u-1/2=\ln(x)$) does not seem to help.

The arguments provided here can be extended to find the density function of
$w=\sqrt{x^{2}/\kappa+\kappa/x^{2}}$ for any $\kappa>0$:%
\[
\frac{w}{\sqrt{2\pi}\sigma}\frac{\exp\left\{  -\dfrac{1}{8\sigma^{2}}\left[
\ln\left(  \dfrac{\kappa}{2}\left(  w^{2}-\sqrt{w^{4}-4}\right)  \right)
-2\mu\right]  ^{2}\right\}  +\exp\left\{  -\dfrac{1}{8\sigma^{2}}\left[
\ln\left(  \dfrac{\kappa}{2}\left(  w^{2}+\sqrt{w^{4}-4}\right)  \right)
-2\mu\right]  ^{2}\right\}  }{\sqrt{w^{4}-4}}
\]
when $w>\sqrt{2}$. \ The case $\kappa=\sqrt{3}$, $\mu=-1/2$ and $\sigma=1$
gives the density function for $a$ in Section 2. \ As another example, if
$\kappa=1$ instead, then%
\[
\operatorname*{E}\left(  w\right)  =3.3278221244164268180344110...,
\]%
\begin{align*}
\operatorname*{E}\left(  w^{2}\right)   &  =\left(  1+e^{2}\right)  e\\
&  =22.8038187516467129762888171....
\end{align*}
Again, an attractive integral representation%
\[
\operatorname*{E}\left(  w\right)  =\dfrac{1}{\sqrt{\pi}}%
{\displaystyle\int\limits_{-\infty}^{\infty}}
\sqrt{\cosh\left(  2v-1\right)  }\exp\left(  -\frac{1}{2}v^{2}\right)  dv
\]
does not seem to help. \ For reasons of brevity, we omit proof of the general formula.

Much more relevant material can be found at \cite{Fi3-LogNrm}, including
experimental computer runs that aided theoretical discussion here.

\section{Addendum}

Continuing a thought raised at the end of Section 3, define a new coordinate
system $\tilde{a}$, $\tilde{b}$, $\tilde{c}$ via a $45^{\circ}$-rotation:%
\[%
\begin{array}
[c]{ccccc}%
a=\dfrac{1}{\sqrt{2}}\left(  \tilde{a}-\tilde{c}\right)  , &  & b=\tilde{b}, &
& c=\dfrac{1}{\sqrt{2}}\left(  \tilde{a}+\tilde{c}\right)
\end{array}
\]
then the orthogonal projection of $\Sigma$ into the plane $\tilde{c}=0$ is the
region bounded (away from the origin) by
\[
\left(  2\tilde{a}^{2}-\tilde{b}^{2}\right)  \tilde{b}^{2}=16.
\]
We note that the boundary is well-approximated, for large $\tilde{a}$, by
$\tilde{b}=\sqrt{2}\tilde{a}$ from above and $\tilde{b}=2\sqrt{2}/\tilde{a}$
from below. \ Ideas for\ \textit{natural} random sampling of $\left(
2\tilde{a}^{2}-\tilde{b}^{2}\right)  \tilde{b}^{2}\geq16$, akin to $a\,b\geq2$
earlier, would be welcome.

\end{document}